\documentclass[10pt]{amsart}
\usepackage{amsfonts,amssymb}
\usepackage[all]{xy}

\newtheorem{thm}{Theorem}[section] 
\newtheorem{lem}[thm]{Lemma} 
\newtheorem{prop}[thm]{Proposition}
\newtheorem{cor}[thm]{Corollary}

\theoremstyle{remark}
\newtheorem{exam}[thm]{Example}

\def\qq{\mathbb{Q}}

\def\rr{\mathbb{R}}
\def\zz{\mathbb{Z}}

\def\cc{\mathbb{C}}
\def\pp{\mathbb{P}}

\def\d{\mathrm{d}}

\def\div{\mathrm{div} \, }

\def\ord{\mathrm{ord}  }
\def\Spec{\mathrm{Spec} \, }

\newcommand*\isom{%
  \xrightarrow{\sim}%
}

\begin{document}

\author{Robin de Jong} 

\title{Local invariants attached to Weierstrass points}

\subjclass[2000]{Primary 14H10, secondary 14H55, 14G40}
\keywords{Faltings height, hyperelliptic
curves, semistable curves, Weierstrass points.}

\maketitle

\begin{abstract} Let $X/S$ be a hyperelliptic curve of genus $g$ over the
spectrum of a discrete valuation ring. Two fundamental numerical invariants are
attached to $X/S$: the valuation $d $ of the hyperelliptic discriminant of
$X/S$,  and the valuation $\delta $ of the Mumford discriminant of $X/S$
(equivalently, the Artin conductor). For a residue field of characteristic~$0$
as well as for $X/S$ semistable the invariants $d $ and $\delta $ are known to
satisfy certain inequalities. We prove an exact formula relating $d $ and
$\delta $ with intersection theoretic data determined by the distribution of
Weierstrass points over the special fiber, in the semistable case. We also prove
an exact formula for the stable Faltings height of an arbitrary curve over a
number field,  involving local contributions associated to its Weierstrass
points. 
\end{abstract}

\section{Introduction} \label{intro}

\thispagestyle{empty}

Let $R$ be a discrete valuation ring with perfect residue field $k$. Let $C$ be
a smooth proper geometrically connected curve of genus $g \geq 1$ over $K=Fr(R)$
and let $\rho \colon X \to S$ be the minimal regular model of $C$ over $S=\Spec
R$. Let $\omega_{X/S}$ be the relative dualising sheaf of $X/S$. The Mumford
isomorphism \cite[Theorem 5.10]{mu}: 
\[ \det R\rho_*(\omega_{X/S}^{\otimes 2})
\otimes K \isom (\det R\rho_* \omega_{X/S} )^{\otimes 13} \otimes K \, , \]
which is well-defined up to a sign, gives a canonical rational section $\Delta$
of: 
\[ (\det R\rho_*\omega_{X/S} )^{\otimes 13} \otimes  \det
R\rho_*(\omega_{X/S}^{\otimes 2})^{\otimes -1} \, , \] 
and a basic numerical
invariant attached to $X/S$, namely the valuation $\delta_{X/S}=\ord_v(\Delta)$
of $\Delta$. Here and in the rest of this paper we fix the convention that the
valuation of a discrete valuation ring is normalised in the sense that its value
group is $\zz$. It turns out that $\delta_{X/S}$ is non-negative and that
actually $\delta_{X/S}=-\mathrm{Art}_{X/S}$, where $\mathrm{Art}_{X/S}$ is the
Artin conductor of $X/S$ \cite[Theorem 1]{sa}. If $X / S$ is semistable one can
view $\Delta$ as the pullback, under the period map $S \to
\overline{\mathcal{M}_g}$, of the tautological section corresponding to the
boundary divisor of the algebraic stack  $\overline{\mathcal{M}_g}$ of stable
curves of genus $g$. In particular if $X/S$ is semistable $\delta_{X/S}$ is 
equal to the number of singular points in the geometric special fiber of $X/S$.

If $C$ is elliptic or hyperelliptic, one has a second basic invariant attached
to $X/S$, namely $d_{X/S}=\ord_v(\Lambda)$ where $\Lambda$ in $H^0(S,(\det
R\rho_* \omega_{X/S})^{\otimes 8g+4})$ is the hyperelliptic discriminant
associated to a hyperelliptic equation for $C/K$  (cf. end of section
\ref{prelim}). Both $\delta_{X/S}$ and $d_{X/S}$ are zero if $X/S$ is smooth.
For $C$ elliptic a result of Tate-Ogg \cite{ogg67} states the equality
$d_{X/S}=\delta_{X/S}$ of both invariants. It follows from examples in 
\cite{ogg66} that in a number of cases with $g=2$ and with $k$ not of
characteristic $2,3$ or $5$, one has an equality $d_{X/S}=2\delta_{X/S}$.
However, as is shown in \cite{ue} there are exceptions to this rule, indicating
that one should add a non-negative `error term' $\varepsilon$ so that the
formula $d_{X/S} = 2 \delta_{X/S} + \varepsilon$ is true in general.  Ueno in
fact compiled a table containing values of $d_{X/S}$ and $\delta_{X/S}$ for $k$
not of characteristic $2,3$ or $5$ based on the complete classification of
degenerate fiber types \cite{nu} in genus two, showing precisely what the
exceptions are. 

A more conceptual approach to the comparison of $d_{X/S}$ and
$\delta_{X/S}$ in genus two was given subsequently by Saito \cite{sa89} and Liu
\cite{liu}. We state the result of Liu. Assume that $k$ is algebraically closed,
and let $\pi \colon X \to X'$ be the contraction of those irreducible components
$\Gamma$ of the special fiber of $X/S$ for which $\deg \omega_{X/S}|_\Gamma=0$.
The hyperelliptic involution $\sigma$ of $X/S$ extends uniquely as an
automorphism of $X'$ over $S$. Let $Z=X'/\langle \sigma \rangle$; this is a
normal surface over $S$, with generic fiber isomorphic to $\pp^1_K$. The
following then holds: let $\tilde{Z} \to Z$ be the minimal desingularisation of
$Z$. Let $n$ be the number of irreducible components of the special fiber of
$\tilde{Z}$. Then $n$ is odd and the formula $ d_{X/S} = 2 \delta_{X/S} + n -1 $
holds. For arbitrary hyperelliptic $C$, Matsusaka \cite[Theorem 4.0.4]{mats} has shown 
the inequalities $ g \delta_{X/S} \leq d_{X/S} \leq g^2 \delta_{X/S} $ 
for $g$ even and $ g \delta_{X/S} \leq d_{X/S} \leq (g^2+1) \delta_{X/S} $ 
for $g$ odd, in the case where $k$ has characteristic~$0$. In particular 
$d_{X/S}$ is non-negative. These inequalities still hold for arbitrary $k$, 
provided one restricts to semistable hyperelliptic curves 
(for definitions see Section \ref{prelim}). In this case the inequalities 
follow from a well-known identity
due, in increasing order of generality,  to Cornalba-Harris \cite{ch}, Kausz
\cite{ka} and  Maugeais \cite{ma} and Yamaki \cite{ya}. 

Our aim in this paper is to prove an exact formula relating $d_{X/S}$ and
$\delta_{X/S}$ in arbitrary genus to intersection theoretic data associated to
the distribution of Weierstrass points over the special fiber of $X/S$. At
present we only have a result for semistable hyperelliptic curves, a
circumstance due the point that we have to assume that all Weierstrass points on
the generic fiber are rational. We hope that this restriction is in the end 
only of a technical nature and can be circumvented by methods as in e.g.
\cite[Section 3]{vi}.

Here is our set-up. Let $F$ be the closure  in $X$  of the Weierstrass points on
$C$. For $P$ a section of $X/S$ we let $\Phi_P$ be the vertical $\qq$-Cartier
divisor on $X$ uniquely determined by the following conditions: the divisor
$(2g-2)P - \omega_{X/S}+\Phi_P$ has intersection number~$0$ with each component
$\Gamma$ of the geometric special fiber of $X/S$, and $P^* \Phi_P$ is the
trivial divisor on $S$. It is known (cf. \cite[Section 9.1.3]{liubook}) that one
has an intersection product $\mathrm{Div}_s X \times \mathrm{Div} X \to \zz$
where $\mathrm{Div} X$ is the group of Cartier divisors on $X$, and
$\mathrm{Div}_s X \subset \mathrm{Div} X$ is the subgroup of Cartier divisors
with support in the special fiber. In particular, one has the intersection
number $(D,\omega_{X/S})$ for $D$ supported in the special fiber. Also, one has
the self-intersection $\Phi_P^2$ of $\Phi_P$ in $\qq$. It is non-positive
\cite[Theorem 9.1.23]{liubook}. 
\begin{thm}  \label{main} Let $R$ be a discrete
valuation ring with perfect residue field and with~$2$ non-zero in $R$. Let $C$
be a hyperelliptic curve of genus $g \geq 2$ with semistable reduction  over
$K=Fr(R)$ and let $X/S$ be the minimal regular model of $C$ over $S=\Spec R$. 
Assume that all Weierstrass points of $C$ are rational over $K$. Then the
integers $d_{X/S}$ and $\delta_{X/S}$ satisfy the following relation: \[ (3g-1)\,
d_{X/S} = - \frac{1}{2} \sum_{P \in F(R)} \Phi_P^2 +  (2g-1)(g+1) \, \delta_{X/S} +
4 \, (E,\omega_{X/S}) \, . \] Here $E$ is the vertical part on $X$ of the divisor of the
Wronskian on an $R$-basis of $H^0(X,\omega_{X/S})$.     
\end{thm}  
The divisor
$E$ will be explained in more detail in Section \ref{weierstrass}. In Section
\ref{calculateE} we will give an explicit formula for $E$. The divisor $\Phi_P$
can be effectively calculated by solving a system of linear equations with
rational coefficients. The assumption on the rationality of all Weierstrass
points already almost guarantees that $C$ has semistable reduction over $K$ (cf.
for a precise statement the  beginning of Section~4 of \cite{ka} and Section
\ref{calculateE} below).  Our theorem has the following global counterpart.
\begin{cor} Let $S$ be a smooth proper connected curve over an algebraically
closed  field $k$ such that $2$ is invertible in $k$. Let $\rho \colon X \to S$
be a semistable hyperelliptic curve of genus $g \geq 2$ with $X$ regular. Assume
that all Weierstrass points of the generic fiber of $X/S$ are rational over the
function field of $S$. Then: 
\begin{align*} (3g-1)(8g+4) & \deg \det R\rho_* \omega_{X/S} \\ & = -
\frac{1}{2} \sum_{P \in F(S)} \Phi_P^2 + (2g-1)(g+1) \sum_{\wp \in |S|}
\delta_\wp + 4 (E, \omega_{X/S}) \, , 
\end{align*}  
where for each closed point $\wp$ of
$S$, $\delta_\wp$ denotes the number of singular points in the fiber of $X/S$ at
$\wp$. The divisor $E$ is defined by taking locally over $S$ the vertical part
of the divisor of the Wronskian on a basis of the locally free $O_S$-module
$\rho_* \omega_{X/S}$.  The divisor $\Phi_P$ is vertical as well and is defined
as in the local case. \end{cor} The proof of Theorem \ref{main} will consist of
a delicate sequence of specialisation and generalisation arguments. The basic
point will be a comparison of the line bundles associated to $d$ and $\delta$ on
the moduli stack of stable hyperelliptic curves. 

A part of what we have to say about Wronskians and Weierstrass points in
families carries over to the more general case of an arbitrary semistable curve
over, say, a Dedekind scheme. In Section \ref{arithmetic} we discuss this point
in an Arakelov context. We obtain a remarkable closed formula for the stable 
Faltings height of a curve over a number field, as well as a lower bound for the
self-intersection of its relative dualising sheaf.  Such lower bounds are of
interest in view of an effective version of the  Bogomolov conjecture
\cite{sz90}. 

In this paper, all schemes are assumed locally noetherian and all morphisms of
schemes are of finite type. A morphism $\rho \colon X \to S$ is called a
prestable curve of genus $g$ if $\rho$ is proper and flat and the geometric
fibers of $\rho$ are connected, reduced, nodal curves of arithmetic genus $g$.
We call a prestable curve semistable (resp. stable) if every smooth rational
component of a geometric fiber meets the other components of that fiber in at
least two (resp. three) points.

\section{Weierstrass points} \label{weierstrass}

We begin by recalling from \cite[Section 3]{ar71} and \cite[Section 2]{vi}  a
basic construction  related to Weierstrass points in families. Let $X,S$ be
regular locally noetherian  schemes and let $\rho \colon X \to S$ be a 
prestable curve of genus $g \geq 1$.  Let $\omega=\omega_{X/S}$ be the relative
dualising sheaf of $X/S$ and let $\lambda$ be the determinant line bundle $\det
\rho_* \omega$ on $S$. Then we claim that the line bundle $\omega^{\otimes
g(g+1)/2} \otimes \rho^* \lambda^{-1}$ on $X$  has a canonical global section
$Wr$. We proceed as follows. Let $X^{sm}$ be the open subset of $X$ where $\rho$
is smooth. Let $x$ be a closed  point on $X^{sm}$ and suppose it maps to a
closed point $s$ of $S$. Let $t$ be a fiber coordinate on an open neighbourhood
$U$ of $x$ in $X^{sm}$. Denote by $\partial^i$ for $i \geq 0$ the
$\hat{O}_{S,s}$-linear selfmap of $\hat{O}_{X,x}=\hat{O}_{S,s}[[t]]$ given by
sending $t^n$ to ${n \choose i} t^{n-i}$. If $O_{S,s}$ is of characteristic~$0$
then $\partial^i$ is just the map sending $f$ to $\frac{1}{i!}\frac{d^i
f}{dt^i}$. For $(f_1,\ldots,f_g)$ in $O_{X,x}^g$ we let
$ Wr(f_1,\ldots,f_g) =
\det \left( \partial^{i-1} f_j \right)_{1 \leq i,j \leq g} $
be the
Wronskian of $(f_1,\ldots,f_g)$. Now let $(\eta_1,\ldots,\eta_g)$ be a basis of
$\rho_* \omega$ around $s$ on $S$ (recall that $\rho_* \omega$ is locally free
of rank $g$), and write $\eta_i = f_i \cdot dt$ with $f_i$ in $O_{X,x}$. If we
then put $ Wr( \eta_1,\ldots,\eta_g ) = Wr(f_1, \ldots, f_g) \cdot (dt)^{\otimes
g(g+1)/2 } $ it can be checked that $Wr(\eta_1,\ldots,\eta_g)$ is independent of
the chosen parameter $t$ and defines a section of  $\omega^{\otimes g(g+1)/2}$
locally around $x$. If $(\eta'_1,\ldots,\eta'_g)$ is another basis of $\rho_*
\omega$ around $s$ connected to $(\eta_1,\ldots,\eta_g)$ via an invertible
linear transformation $(a_{ij})$  then one easily verifies that
$Wr(\eta'_1,\ldots,\eta'_g) = \det(a_{ij}) Wr(\eta_1,\ldots,\eta_g)$. This
implies that if we take $ Wr(\eta_1,\ldots,\eta_g) \cdot (\eta_1 \wedge \ldots
\wedge \eta_g)^{-1} $ we obtain a local section of $\omega^{\otimes g(g+1)/2}
\otimes \rho^* \lambda^{-1}$ that is independent of the choice of a basis. As we
can glue these local sections over $X^{sm}$, and extend uniquely over $X$, we
obtain a canonical section $Wr$ of $\omega^{\otimes g(g+1)/2} \otimes \rho^*
\lambda^{-1}$ as required. 

We call $Wr$ the Wronskian of $X/S$.  For $X/S$
smooth the construction of the Wronskian commutes with \'etale base change. If
$S = \Spec k$ where $k$ is a field and $Wr$ is not identically zero we call the
divisor $\div Wr$ of $Wr$ the divisor of Weierstrass points of $X/S$.  In
general, if the Wronskian is not identically zero, we call the closure in $X$ of
the Weierstrass points of the general fibers of $X/S$  the Weierstrass divisor
of $X/S$ and denote it by $W$. We can then write $\div Wr = W+E$ for some
effective divisor $E$ on $X$ where we call $E$ the residual divisor of $X/S$.
Both $W$ and $E$ can be viewed as Cartier divisors on $X$. 

The next theorem
slightly generalises Lemma 3.3 of \cite{ar71}  and Theorem 2.10 of \cite{vi}. 
\begin{prop} \label{canonicalisom} (Arakelov) 
Let $\rho \colon X \to S$ be a prestable
curve of genus $g \geq 1$ with $X,S$ regular locally noetherian schemes.  Let
$Wr$ be the Wronskian of $X/S$.  If $Wr$ is not identically zero there exists a
canonical isomorphism:  
\[ \omega_{X/S}^{\otimes g(g+1)/2} \otimes \rho^* (\det
\rho_* \omega_{X/S})^{-1} \isom O_X(W+E)  \] 
of line bundles on $X$  with $W$
the Weierstrass divisor of $X/S$ and with $E$ the residual divisor of $X / S$.
$E$ has support in the union of the non-smooth fibers of $X/S$ and the fibers 
of positive characteristic $p$ with $p < 2g-1$. For $X/S$ smooth the  formation
of the isomorphism commutes with any \'etale surjective  base change. 
\end{prop}  
\begin{proof} All statements are immediate from our discussion,
except perhaps for the point that for $X/S$ smooth, the divisor $E$ may have
support in the fibers of positive characteristic $p$ with $p < 2g-1$. For this
we may assume that $S=\Spec k$ with $k$ an algebraically closed field. Let $p$
be the characteristic of $k$. We would like to prove that if $Wr$ is zero, then
$p$ is positive with $p < 2g-1$. Let $x$ be a closed point of $X$, and let $G(x)
= \{ a_1 < a_2 < ... \}$ be the set of natural numbers $a$ such that there is an
$\eta$ in $H^0(X,\omega)$ with a zero of exact order $a-1$ at $x$. Then $G(x)$
gives rise to a filtration with simple quotients of the $g$-dimensional
$k$-vector space $H^0(X,\omega)$, so that $G(x)$ consists of $g$ elements. Since
$\omega$ has no base points we have $a_1=1$. On the other hand we have $a_g \leq
2g-1$. Moreover there is a local parameter $t$ around $x$ and a basis
$(f_1dt,\ldots,f_gdt)$ of $H^0(X,\omega)$ such that  $f_i = (\mathrm{unit})
\cdot   t^{a_i-1} $ for $i=1,\ldots,g$. This implies that there is an element
$A$ in the image of $\zz$ in $O_{X,x}$ such that: 
\begin{align*}  A \cdot
Wr(f_1,\ldots,f_g) & =  \left| \begin{array}{cccc}  t^{a_1-1} &
(a_1-1)t^{a_1-2} & \ldots & (a_1-1)\cdots(a_1-g+1)t^{a_1-g} \\  t^{a_2-1} &
(a_2-1)t^{a_2-2} & \ldots & (a_2-1)\cdots(a_2-g+1)t^{a_2-g} \\ \vdots    &
\vdots           &        & \vdots   \\ t^{a_g-1} & (a_g-1)t^{a_g-2} & \ldots &
(a_g-1)\cdots(a_g-g+1)t^{a_g-g} \end{array} \right| \cdot (\mathrm{unit})   \\ &
=  \left| \begin{array}{cccc}  1 & (a_1-1)  & \ldots & (a_1-1)\cdots(a_1-g+1) 
\\  1 & (a_2-1)  & \ldots & (a_2-1)\cdots(a_2-g+1)  \\ \vdots    &
\vdots           &        & \vdots   \\ 1 & (a_g-1)  & \ldots &
(a_g-1)\cdots(a_g-g+1)  \end{array} \right| \cdot (\mathrm{unit}) \cdot t^w   \\
& =  \prod_{1 \leq i < j \leq g} (a_i -a_j)  \cdot (\mathrm{unit})  \cdot t^w
\end{align*} 
where $w = \sum_{i=1}^g (a_i-i)$. Thus if
$Wr(f_1dt,\ldots,f_gdt)$ is zero, so is $\prod_{1 \leq i < j \leq g} (a_i
-a_j)$, whence  $k$ of positive characteristic $p$ with $p<2g-1$. \end{proof}
The divisor $W$ is in general not finite over $S$, nor is the decomposition of
$\div Wr$ as $W+E$ in general well-behaved under base change (cf. \cite[Remark
2.11]{vi}). The situation is better if $\dim S=1$ or if $X/S$ is a semistable
hyperelliptic curve (for which see the next section).

\section{Hyperelliptic curves}

\subsection{Preliminaries} \label{prelim}

We begin by defining precisely what we mean by a smooth hyperelliptic curve over
a base scheme $S$. A detailed account can be found in \cite{lk}. Let $k$ be an
algebraically closed field, and let $X/k$ be a proper smooth curve. We call
$X/k$ hyperelliptic if the genus $g$ of $X$ is at least two and there is an
involution $\sigma$ in $\mathrm{Aut}_k(X)$ such that $X/\langle \sigma \rangle
\cong \pp^1$. Such an automorphism is then unique. If $S$ is a locally
noetherian scheme and $\rho \colon X \to S$ is a smooth proper curve we say that
$X/S$ is hyperelliptic if there exists an involution $\sigma$ in
$\mathrm{Aut}_S(X)$ such that the restriction of $\sigma$ to each geometric
fiber of $\rho$ gives that fiber  the structure of a hyperelliptic curve (cf.
\cite[Definition 5.4 and  Theorem 5.5]{lk}). We note that the involution
$\sigma$ is uniquely determined since the automorphism scheme
$\underline{\mathrm{Aut}}_S(X)$ is unramified over $S$ (\cite[Theorem
1.11]{dm}). We call $\sigma$ the hyperelliptic involution of $X/S$.  If $X/S$ is
a smooth hyperelliptic curve, we denote by $F=F_{X/S}$ the fixed point subscheme
of $X$ under the action of the group $\langle \sigma \rangle$. By
\cite[Proposition 6.3]{lk} the scheme $F$ is the closed  subscheme associated to
an effective Cartier divisor on $X$ relative to $S$. The divisor $F$ is finite
and flat over $S$ of degree $2g+2$,  and its formation commutes with arbitrary
base change. 
\begin{lem} \label{smoothhyp} Assume that $S$ is a regular locally
noetherian scheme,  and that $X/S$ is a smooth hyperelliptic curve of genus $g$.
Then $X$ is regular, the Wronskian $Wr$ of $X/S$ is not identically zero, the
residual divisor $E$ of $X/S$ is empty, and the equality of Cartier divisors $W=
g(g-1)/2 \cdot F$ holds. 
\end{lem}       
\begin{proof} That $X$ is regular
follows from \cite[Theorem 4.3.36]{liubook}. For the rest it suffices to
consider the case $S =\Spec k$ with $k$ an algebraically closed field. Then an
affine part of  $X/S$ can be given by an equation $y^2+ay=b$ with $a,b$ in
$k[x]$. The hyperelliptic involution $\sigma$ is given by $y \mapsto -y-a$. As
$X/\langle \sigma \rangle$ has genus~$0$, the quotient map $X \to X/\langle
\sigma \rangle$ is separable. This implies that $2y+a$ is not identically zero,
and that a basis of the regular differentials is given by $x^i \d x/(2y+a)$ for
$i=0,\ldots,g-1$.  A computation yields: \[ Wr \left( \frac{\d x}{2y+a}, \ldots,
\frac{ x^{g-1} \d x}{2y+a} \right) = \left( 2y+a \right)^{g(g-1)/2} \left(
\frac{\d x}{2y+a} \right)^{\otimes g(g+1)/2} \, , \] which is not identically
zero. It follows that $E$ is empty and finally, since $F$ is locally given by
the vanishing of $2y+a$, the equality $W=g(g-1)/2 \cdot F$. 
\end{proof}  
If $\rho \colon X \to S$ is an arbitrary proper flat curve with $S$ locally
noetherian, we call $X/S$ hyperelliptic if there exists an involution $\sigma$
in $\mathrm{Aut}_S(X)$ and an open dense subset $U$ of $S$ such that the
restriction of $\rho$ to $X \times_S U$ is a smooth hyperelliptic curve with the
restriction of $\sigma$ to $X \times_S U$ as hyperelliptic involution. We remark
that if $S$ is a connected regular noetherian scheme of dimension~$1$ and  $X$
is the regular minimal or stable model of its generic fiber, then $X/S$ is
hyperelliptic if and only if its generic fiber is hyperelliptic. Indeed, a
hyperelliptic involution of the generic fiber extends by the valuative criterion
of properness, and according to \cite[Proposition 5.14]{lk} the smooth fibers of
$\rho$ are then organised in a smooth hyperelliptic curve over an open dense
subset of $S$. We denote by $F$ the closure in $X$ of the fixed point subscheme 
on the smooth fibers of $X/S$.

Let $S$ again be arbitrary.  Let $\rho \colon X \to S$ be hyperelliptic of genus
$g \geq 2$,  and assume $X$ to be regular. Let $\omega_{X/S}$ be the relative
dualising sheaf of $X/S$. The line bundle $\left( \det \rho_* \omega_{X/S}
\right)^{\otimes 8g+4}$ on $S$ has then a canonical non-zero global section
$\Lambda$,  called the hyperelliptic discriminant of $X/S$ \cite[Section
2]{ka}.  For $S = \Spec R$ with $R$ a discrete valuation ring with $2$ non-zero
in $R$ it can be defined as follows. Let $K=Fr(R)$, and let $y^2=f(x)$ with $f$
in $K[x]$ monic separable of degree $2g+2$  be an equation for $X_K$.  Then we
put: 
\[ \Lambda = \left( 2^{-(4g+4)} \cdot D(f) \right)^g \cdot \left( \frac{\d
x}{y} \wedge \ldots \wedge \frac{x^{g-1} \d x}{y} \right)^{\otimes 8g+4} \, , \]
where $D(f)$ is the discriminant of $f$ in $K$. The section $\Lambda$  is
independent of the choice of $f$ and hence defines a non-zero canonical element
of $H^0(X_K,\Omega^1_{X_K})^{\otimes 8g+4}$. It can then be viewed as a rational
section of $\left( \det \rho_* \omega_{X/S} \right)^{\otimes 8g+4}$ over $S$.
The integer $d_{X/S}$ is  defined to be the valuation $\ord_v(\Lambda)$ of
$\Lambda$ at the closed point $v$ of $S$. The construction of $\Lambda$
generalises over arbitrary locally noetherian base schemes \cite[Proposition
2.7]{ma}. For $X/S$ smooth the section $\Lambda$ is nowhere vanishing on $S$.

\subsection{Proof of Theorem \ref{main}}

The proof of Theorem \ref{main} goes by a number of  specialisation and
generalisation arguments. Whenever $\rho \colon X \to S$  is a semistable curve,
we just write $\omega$ for its relative dualising sheaf $\omega_{X/S}$ and
$\lambda$ for its determinant line bundle $\det R\rho_* \omega$. The formation
of both $\omega$ and $\lambda$ commutes with arbitrary base change. When dealing
with the tensor product of line bundles, we will often mix additive and
multiplicative notation; we hope that this does not lead to confusion.

We start with an arbitrary generically smooth semistable curve $\rho \colon X
\to S$ of genus $g \geq 2$ with $S$ locally noetherian. We consider the line
bundle \[ Q= \langle (2g-1)(g+1)\omega - 4 \rho^* \lambda, \omega \rangle \] on $S$;
here $\langle L, M \rangle$ denotes the Deligne symbol of two line bundles $L,M$
on $X$ \cite[Section 6 and 7]{de} \cite[Section 5.4]{mbell}. The Deligne symbol
is a biadditive functorial pairing from the category of line bundles on $X$ to
the category of line bundles on $S$. For $L=\rho^*N$ and for $M$ flat over $S$
one has $\langle L,M \rangle = \langle \rho^* N,M \rangle \isom N^{\otimes d}$
canonically, where $d$ is the degree of $M$ in the fibers of $\rho$. In our case
we obtain a canonical isomorphism of line bundles $Q \isom (2g-1)(g+1) \langle
\omega, \omega \rangle - (8g-8)\lambda$ on $S$. According to \cite[Th\'eor\`eme
2.1]{mb2} we also have a canonical isomorphism of line bundles $\langle
\omega,\omega \rangle \isom 12 \lambda - \delta$ on $S$, where $\delta$ is the
line bundle on $S$ given by the boundary divisor on the moduli stack of stable
curves of genus $g$ or, equivalently, the zero divisor of the section $\Delta$
coming from the Mumford isomorphism over the locus of $S$ where $\rho$ is
smooth.  A combination of these two isomorphisms yields a canonical isomorphism
\[ (*) \qquad Q \isom (3g-1)(8g+4)\lambda -(2g-1)(g+1)\delta \] of line bundles on $S$. Note that all isomorphisms are compatible with any base change preserving the generic
smoothness of $\rho$.

Next we impose on $\rho$  the condition that it is hyperelliptic, and on $S$ the
condition that it is both reduced and irreducible, and that it has a non-empty
open subscheme $U$ which is both regular and a scheme over $\zz[1/2]$. We claim
that under these assumptions the line bundle $Q^{\otimes 2}$ has a canonical
non-zero rational section $\xi$, having the property that its valuation
$\ord_v(\xi)$ at the closed point $v$ of the spectrum $S$  of a discrete
valuation ring equals $\ord_v(\xi)= -\sum_{P \in F(S)} \Phi_P^2 + 8 (E,
\omega)$, if all Weierstrass points on the generic fiber of $X/S$ are rational,
and $X$ is the regular minimal model of its generic fiber. As $S$ is reduced and
irreducible, it suffices to construct $\xi$ over the non-empty open subscheme 
$U$. By shrinking $U$ if necessary, we may assume that $\rho$ is smooth over
$U$. 

Hence we reduce to the case that $\rho$ is smooth over $S$ and that $S$ is
regular and a scheme over $\zz[1/2]$. In this case $X$ itself is regular by
Lemma \ref{smoothhyp} and the fixed point subscheme $F$ of $X$ over $S$ is
finite \'etale of degree $2g+2$ over $S$  by \cite[Corollary 6.8]{lk}. After a
faithfully flat base change we can assume that $F$ is the disjoint union of
$2g+2$ sections of $S$. By   faithfully flat descent we may  assume that this is
then already the case. Thus, let $P$ be a section of $F$ over $S$. By
\cite[Lemma 6.2]{dJ07} we have a unique isomorphism $\omega \isom (2g-2)O_X(P)
\otimes \rho^* \langle P,P \rangle^{\otimes -(2g-1)}$ of line bundles on $X$ 
that induces, by pulling back along $P$, the canonical adjunction isomorphism
$\langle P, \omega \rangle \isom \langle P,P \rangle ^{\otimes -1}$ (cf.
\cite[Th\'eor\`eme 6.11]{mb85}).  

The formation of this isomorphism commutes
with arbitrary base change. Pairing the isomorphism with $\omega$ and using
adjunction one obtains a canonical isomorphism $\langle \omega, \omega \rangle
\isom 4g(g-1) \langle P, \omega \rangle$ of line bundles on $S$. Taking the sum
over all $P$ we obtain, using that $W = g(g-1)/2 \cdot F$ by Lemma
\ref{smoothhyp}, a canonical isomorphism $\langle 4W - (g+1)\omega,\omega
\rangle^{\otimes 2} \isom O_S$ of line bundles on $S$. With Lemma
\ref{canonicalisom} we get, as $E$ is empty in the present case by Lemma
\ref{smoothhyp}, a canonical isomorphism $Q^{\otimes 2}= \langle
(2g-1)(g+1)\omega - 4\rho^* \lambda , \omega \rangle^{\otimes 2} \isom O_S$. We
get $\xi$ by taking the section of $Q^{\otimes 2}$ that corresponds to the
canonical section~$1$ of $O_S$ under this isomorphism. 

Now assume that $S$ is the spectrum of a discrete valuation ring, that all
Weierstrass points on the generic fiber of $X/S$ are rational, and that $X$ is
the regular minimal model of its generic fiber. It follows from the construction
of $\xi$ above that for each $P\in F(S)$ we have a canonical trivialising 
section $s_P$ of the line bundle $ \omega - (2g-2)O_X(P) + \rho^* \langle P, P
\rangle^{\otimes 2g-1} $ restricted to the generic fiber of $X/S$. This section
extends uniquely to a rational section $s_P$ of this same line bundle over $X$.
Let's denote by $\Phi'_P$ its divisor. Then $\Phi'_P$ is supported on the
special fiber of  $X/S$ and we have $P^* \Phi'_P$ trivial on $S$ since the
pullback of our  line bundle along $P$ is $\langle P, P \otimes \omega \rangle$
which is canonically trivial by adjunction. It follows easily that  $\Phi'_P$ is
equal to  $\Phi_P$. Thus we have a canonical isomorphism \[ \omega - (2g-2)O_X(P)
+ \rho^* \langle P,P \rangle^{\otimes 2g-1 } \isom O_X(\Phi_P)  \] of line
bundles on $X$. Pairing with $\omega$ yields a canonical isomorphism \[ -4g(g-1)
\langle P, \omega \rangle + \langle \omega , \omega \rangle \isom \langle
\Phi_P, \omega \rangle  \] of line bundles on $S$. On the other hand, pairing
with $\Phi_P$ yields the isomorphism  $ \langle \Phi_P, \omega \rangle \isom
\langle \Phi_P, \Phi_P \rangle$. Combining we find an isomorphism  \[  -4g(g-1)
\langle P, \omega \rangle + \langle \omega , \omega \rangle \isom  \langle
\Phi_P, \Phi_P \rangle \] and summing over $P \in F(S)$ and adding $8\langle
E,\omega \rangle$  on both sides gives the isomorphism  $  \langle 4W -
(g+1)\omega + 4E, \omega \rangle^{\otimes 2}  \isom - \sum_{P \in F(S)} \langle
\Phi_P, \Phi_P \rangle  + 8 \langle E, \omega \rangle $. Using Lemma
\ref{canonicalisom} we find a canonical isomorphism \[Q^{\otimes 2} = \langle
(2g-1)(g+1)\omega-4\rho^*\lambda,\omega \rangle^{\otimes 2} \isom - \sum_{P \in
F(S)} \langle \Phi_P, \Phi_P \rangle  + 8 \langle E, \omega \rangle \, . \] The line
bundle on the right hand side has a canonical non-zero rational section $\xi'$
whose valuation at the closed point $v$ of $S$ equals $\ord_v(\xi')= -\sum_{P
\in F(S)} \Phi_P^2 + 8 (E,\omega)$. It follows from the constructions that $\xi$
in $Q^{\otimes 2}$ corresponds to $\xi'$, and thus we find the required formula
for $\xi$. 

The proof of Theorem \ref{main} will be done if we could also prove the formula
$\ord_v(\xi) = (6g-2)\ord_v(\Lambda) - (4g-2)(g+1)\ord_v(\Delta)$ for the
valuation at $v$ of $\xi$. For this, recall the canonical isomorphism $Q \isom
(3g-1)(8g+4)\lambda -(2g-1)(g+1)\delta$ of line bundles on $S$ that we obtained
from Mumford's isomorphism. View $\Lambda^{\otimes 3g-1}$ as a rational section
of the right hand side of this isomorphism. We will be done once we prove that
under the square of the isomorphism (*) the rational section $\xi$ is identified
with $\Lambda^{\otimes 6g-2}$, up to a sign. In order to accomplish this, we
first prove that for each stable hyperelliptic curve $\rho \colon X \to S$ with
$S$ locally noetherian, reduced and irreducible and having the property that it
contains a non-empty open subscheme $U$ that is regular and defined over
$\zz[1/2]$, the sections $\xi$ and $\Lambda^{\otimes 6g-2}$ are identified, up
to a sign. 

For this it suffices to consider the algebraic stack
$\overline{\mathcal{I}_g}$ classifying stable hyperelliptic curves of genus $g$ 
studied in e.g. \cite[Section 4b]{ch} (over $\cc$) and \cite[Section 1]{ya}
(over $\zz$). It is a suitable compactification of the algebraic  stack
$\mathcal{I}_g$ of smooth hyperelliptic curves of genus $g$. The stack
$\mathcal{I}_g$ has smooth and geometrically irreducible fibers over $\Spec \zz$
(cf. \cite[Theorem 3]{ll}) and is in particular itself reduced and irreducible.
The same holds then for $\overline{\mathcal{I}_g}$ and it is clear that the
latter contains a non-empty open substack which is  regular and defined over
$\zz[1/2]$ (just take $\mathcal{I}_g \otimes \zz[1/2]$). For any stable
hyperelliptic curve $\rho \colon X \to S$ of genus $g$ with $S$ locally
noetherian, reduced and irreducible and having a non-empty regular open
subscheme defined over $\zz[1/2]$, we obtain the sections $\xi$ and
$\Lambda^{\otimes 6g-2}$ over $S$  by pullback from $\overline{\mathcal{I}_g}$
under the period map  $S \to \overline{\mathcal{I}_g}$. 

Now note that both $\xi$ and $\Lambda^{\otimes 6g-2}$ are supported on the
boundary of $\mathcal{I}_g$ in $\overline{\mathcal{I}_g}$. Letting $\phi$ be the
square of the canonical isomorphism (*) over $ \overline{\mathcal{I}_g}$ we find
that $\phi(\xi) \otimes \Lambda^{\otimes -(6g-2)}$ is a rational function on
$\overline{\mathcal{I}_g}$, regular invertible  on $\mathcal{I}_g$. By
\cite[Proposition 7.3]{dJ07}, the function  $\phi(\xi) \otimes \Lambda^{\otimes
-(6g-2)}$ equals $\pm 1$. Thus the sections $\xi$ and $\Lambda^{\otimes 6g-2}$
are identified, up to a sign, as required. 

The proof will now be finished by making the transition, in the case that $S$ is
the spectrum of a discrete valuation ring, from a stable model $X/S$ to its
associated regular minimal model $X'/S$ with structure morphism $\rho'$.  Let
$\omega'$ be the relative dualising sheaf of $X'/S$, and write \[ \lambda'=\det
\rho_* \omega' \, , \quad Q'=\langle (2g-1)(g+1)\omega'-4\rho'^*\lambda',\omega'
\rangle \, . \] Let $\xi'$ be the canonical rational section of $Q'^{\otimes 2}$, and
let $\Lambda'$ be the canonical rational section of $\lambda'^{\otimes 8g+4}$.
Finally let $\pi \colon X' \to X$ be the canonical map contracting the
$(-2)$-curves in the special fiber.  We know that for $X/S$ the sections $\xi$
of $Q^{\otimes 2}=\langle (2g-1)(g+1)\omega-4\rho^*\lambda,\omega \rangle$ and
$\Lambda^{\otimes 6g-2}$ of $(6g-2)(8g+4)\lambda - (4g-2)(g+1)\delta$
correspond, up to a sign. We have a canonical isomorphism $\pi^* \omega \isom
\omega'$ and for any pair of line bundles  $L,M$ on $X$ a canonical isomorphism
$\langle \pi^*L, \pi^*M \rangle \isom \langle L,M \rangle$ \cite[Section
5.4]{mbell}.  Thus we find a canonical isomorphism $Q \isom Q'$, and one
verifies that $\xi$ and $\xi'$ are identified in this way.  Also we have a
canonical isomorphism $\lambda \isom \lambda'$, yielding an  identification of
$\Lambda$ and $\Lambda'$. We conclude that $\xi'$ and $\Lambda'^{\otimes 6g-2}$
are identified as well.

\subsection{Effective computations} \label{calculateE}

In this section we indicate how, for a semistable hyperelliptic curve $X$ over a
discrete valuation ring $R$ in which  $2$ is a unit, the residual divisor $E$
can be effectively calculated. We base our discussion on Section~4 of I. Kausz's
article \cite{ka}. Another approach can be found in \cite[Section 4]{vi}. We
assume that all Weierstrass points on the generic fiber are rational. This
implies (cf. \cite[Lemma 4.1]{ka}) that the generic fiber of $X/R$ has an
equation $y^2 = Af(x)$ with $A$ a unit in $R$ and  $f(x)= \prod_{i=1}^{2g+2}
(x-a_i)$ for certain pairwise distinct  $a_i \in R$. At the expense of making a
small finite extension of $R$ we assume that the $v(a_i-a_j)$ are even for all
$i \neq j$, and that the number of distinct images of the $a_i$ in $R/m$ is at
least~$3$. Here $v$ is the normalised discrete valuation of $R$, and $m$ is the
maximal ideal of $R$.

Assume that $X$ is the minimal regular model of its generic fiber.  We can
describe $X/R$ in a combinatorial way. We start by constructing a finite tree
$T=(\mathcal{V},\mathcal{E})$ from the above data.  For non-negative integers
$n$ denote by  $r_n \colon \{a_1,\ldots,a_{2g+2} \} \to R/m^n$ the natural map
sending $a_i$ to its residue class modulo $m^n$. The vertices of $T$ are then
the elements of the set $\mathcal{V} = \sqcup_{n \geq 0} \mathcal{V}_n$ where
$\mathcal{V}_n = \{ V \in R/m^n \, \colon \, \# r_n^{-1}(V) \geq 2 \}$. The set
$\mathcal{E}$ of edges of $T$ consists of the pairs $(V,V')$ where $V \in
\mathcal{V}_n$ and $V' \in \mathcal{V}_{n+1}$ for some $n \geq 0$ and $V'
\mapsto V$ under the canonical map $\mathcal{V}_{n+1} \to \mathcal{V}_n$. It
follows that $\mathcal{V}$ has a canonical partial ordering and that there is a
unique minimal element  $V_0$ with respect to this ordering. Moreover, $T$ is
canonically isomorphic to the dual graph of the special fiber of the prestable
curve $Y' /R$ of genus~$0$ that is obtained by taking the smooth curve $\pp^1_R$
and then successively blowing up the closed points of the special fiber where
the sections $P_i$ given by the $a_i$ meet, until the strict transform of
$\sum_i P_i$ in $Y'$ becomes regular.

We can construct $X$ from $Y'$. For every $V \in \mathcal{V}$ put $n(V)=n$ if $V
\in \mathcal{V}_n$, and put $\varphi(V) = \#r_n^{-1}(V)$.  Next define $C(V)$ to
be~$1$ if both $n$ and $\varphi(V)$ are odd, and~$0$ otherwise.  This gives rise
to an effective divisor $C = \sum_{i=1}^{2g+2} P_i + \sum_{V \in \mathcal{V}}
C(V) \cdot V$ on $Y'$ which has the properties that $C$ is regular and that the
class of $C$ is divisible by two in the Picard group of $Y'$. By standard
constructions we obtain an  $R$-scheme $X'$ and a finite flat morphism $\pi'
\colon X' \to Y'$ of degree two such that $X'$ is regular and $\pi'$ is branched
exactly along $C$. In fact $X' /R$ is prestable with generic fiber isomorphic to
the generic fiber of $X$. For $V$ an irreducible component of the special fiber
of $Y'$, set $\widetilde{V} = \pi'^*V = X' \times_{Y'} V$. If $C(V)=1$ then
$\widetilde{V}=2L$ with $L$ an exceptional smooth rational curve. If $C(V)=0$
then $\widetilde{V}$ is reduced and $\widetilde{V} \to V$ is finite of
degree two, ramified over precisely the $(V,C)$ intersection points of  $C$ with
$V$. Upon contracting all exceptional smooth rational curves in the special
fiber  of $X'$ we find the regular semistable model $X / R$, up to 
$R$-isomorphism. Likewise, one can contract all components $V$ on $Y'$ with
$C(V)=1$. We denote the resulting  model by $Y$, and we have a canonical map
$\pi \colon X \to Y$ of $R$-schemes.  
\begin{prop} \label{formulaE} Let $X$ be a
hyperelliptic semistable curve of genus $g \geq 2$ over a  discrete valuation
ring $R$ in which~$2$ is a unit. Assume that $X$ is regular.  Suppose that the
generic fiber of $X /R$ is given by an equation $y^2 = A \cdot f(x)$ with $A \in
R^*$ and $f(x)= \prod_{i=1}^{2g+2} (x-a_i)$ for  certain distinct $a_i \in R$.
Assume that $v(a_i-a_j)$ is  even for $i \neq j$, and that the number of
distinct images of the $a_i$ in $R/m$ is at least~$3$. Let
$T=(\mathcal{V},\mathcal{E})$ be the tree associated to the $a_i$ as above, and
let $e$ be the integer:   \[  e = \frac{1}{2} \sum_{V>V_0 \colon \varphi(V) \,
\mathrm{even}} \frac{\varphi(V)}{2}\left( \frac{\varphi(V)}{2}-1 \right) +
\frac{1}{2} \sum_{V>V_0 \colon \varphi(V) \, \mathrm{odd}} \left(
\frac{\varphi(V)-1}{2} \right)^2 \,  .\]  Here the sums run over $V$ in the
vertex set $\mathcal{V}$. Then for the residual divisor $E$ of $X/R$ one has the
formula:  \[ E = \sum_{V \in \mathcal{V} \atop C(V) = 0} \left( e - \frac{g}{2}
\sum_{i=1}^{n(V)} \varphi(V_i) + \frac{g(g+1)}{2} n(V) \right) \cdot
\widetilde{V} \, , \] where for each given $V$ in $\mathcal{V}$ we denote by
$V_0,V_1,\ldots,V_n=V$ the vertices of the unique linear subgraph of $T$ that
connects $V$ and $V_0$. 
\end{prop}  
\begin{proof} The divisor $E$ is the
vertical part of the divisor on $X$ of the  Wronskian
$Wr(\omega_0,\ldots,\omega_{g-1})$ on an  $R$-basis
$(\omega_0,\ldots,\omega_{g-1})$ of $H^0(X,\omega)$. According to
\cite[Proposition 5.5]{ka} there are $e_i \in \zz$ with $\sum_{i=0}^{g-1} e_i =
e$ and $b_j \in \{a_1,\ldots,a_{2g+2} \}$  such that  for: \[ \omega_i =
t^{e_i}  \left( \prod_{j=1}^i (x-b_j) \right) \frac{\d x}{y} \qquad , \qquad
i=0,\ldots,g-1\]  the tuple $(\omega_0,\ldots,\omega_{g-1})$ is an $R$-basis of
$H^0(X,\omega)$.  Here $t$ is a generator of the maximal ideal $m$ of $R$.
Putting $h_i(x) = \prod_{j=1}^i(x-b_j)$ for $i=0,\ldots,g-1$ a computation shows
that:   \begin{align*} Wr(\omega_0,\ldots,\omega_{g-1}) & =  y^{-g} t^e
Wr(h_0,\ldots,h_{g-1})  (\d x)^{\otimes g(g+1)/2}  \\ & =  y^{g(g-1)/2} t^e
Wr(h_0,\ldots,h_{g-1}) \left( \frac{\d x}{y} \right)^{\otimes g(g+1)/2} \\ & =
  y^{g(g-1)/2} t^e \left( \frac{\d x}{y} \right)^{\otimes g(g+1)/2} \, . 
\end{align*}  From this we compute $ \div Wr( \omega_0,\ldots,\omega_{g-1}
)$. Let $\mathrm{div}_\mathrm{vert} y$ be the vertical part of the divisor of
$y$ on $X$, and let $P_\infty$ be the section of $Y$ corresponding to the point
at infinity of the generic fiber of $Y$. According to Lemma 5.2 of \cite{ka} we
have:  \[ \div \left( \frac{\d x}{y} \right) = (g-1)\pi^*P_\infty -
\mathrm{div}_\mathrm{vert} y + \sum_{V \in \mathcal{V} \atop C(V)=0} n(V) \cdot
\widetilde{V} \, , \] where $\pi \colon X \to Y$ is the canonical map. It
follows that:  \begin{align*} \div & Wr( \omega_0,\ldots,\omega_{g-1} )  =  e
F + \frac{g(g-1)}{2} \div y + \frac{g(g-1)(g+1)}{2} \pi^*P_\infty \\  &  -
\frac{g(g+1)}{2}  \mathrm{div}_\mathrm{vert} y +  \frac{g(g+1)}{2} \sum_{V \in
\mathcal{V} \atop C(V)=0} n(V) \cdot \widetilde{V}  \end{align*}  with $F$
the special fiber of $X$. Noting that:  \[ \frac{g(g-1)}{2} \div y =
\frac{g(g-1)}{2} \mathrm{div}_\mathrm{vert} y  + W -  \frac{g(g-1)(g+1)}{2}
\pi^*P_\infty \]  we derive:  \[ \div Wr( \omega_0,\ldots,\omega_{g-1} ) = e  F
+ W - g \mathrm{div}_\mathrm{vert} y  + \frac{g(g+1)}{2} \sum_{V \in \mathcal{V}
\atop C(V)=0} n(V) \cdot \widetilde{V} \, . \] Clearly for each $V$ in
$\mathcal{V}$ with $C(V)=0$ we have $ v_{\widetilde{V}}(y) = \frac{1}{2}
\sum_{i=1}^{2g+2} v_{\widetilde{V}} (x-a_i) $. The valuation $v_{\widetilde{V}}
(x-a_i)$ can be seen to be equal to $\min(n(V),v(a-a_i))$ where $a$ is a
representative of $V$ (cf. \cite[Proof of Lemma 5.1]{ka}). A counting argument
shows that  $ \sum_{i=1}^{2g+2} \min(n(V),v(a-a_i))  = \sum_{i=1}^{n(V)}
\varphi(V_i) $ if $V_0,V_1,\ldots,V_n=V$ are the vertices of the unique linear
subgraph of $T$ that connects $V$ and $V_0$. Thus we obtain 
$v_{\widetilde{V}}(y) = \sum_{i=1}^{n(V)} \varphi(V_i) $ for each $V$ with
$C(V)=0$ and we arrive at:  \[  \div Wr( \omega_0,\ldots,\omega_{g-1} ) = W +
\sum_{V \in \mathcal{V} \atop C(V)=0}   \left( e - \frac{g}{2} \sum_{i=1}^{n(V)}
\varphi(V_i) + \frac{g(g+1)}{2} n(V) \right) \cdot \widetilde{V} \, . \]  The
formula follows.  
\end{proof}  
Apparently, the number $e$ can be interpreted as
the multiplicity in $E$  of the irreducible component of the special fiber of $X
/ R$ that maps  to $V_0$ in $Y/R$. We note however that the component $V_0$ may
depend on the particular equation chosen for the generic fiber of $X / R$.

\subsection{Example}

In this section we verify Theorem \ref{main} for a concrete case with $g=2$.
Let $k$ be an algebraically closed field with $2 \in k^*$, and let $X'_k$
be a stable curve of genus two over $k$ consisting of an elliptic curve $A$ with
a one-noded  rational curve $B$ attached to it. Let $X_k \to X'_k$ be the
modification of $X'_k$ obtained by partially normalising $X'_k$ at the node
$\nu$ of $B$, and attaching a projective line $D$ at the two points in the 
preimage of $\nu$. The semi-stable curve obtained in this way has type
$I_2-I_0-1$ in the classification of Namikawa-Ueno \cite{nu}. Its intersection
matrix is: \\
\smallskip
\medskip\noindent
\vbox{
\bigskip\centerline{\def\quad{\hskip 0.6em\relax}
\def\quod{\hskip 0.5em\relax }
\vbox{\offinterlineskip
\hrule
\halign{&\vrule#&\strut\quod\hfil#\quad\cr
height2pt&\omit&&\omit&&\omit&&\omit&\cr
&    &&  $A$    &&   $B$ && $D$       &\cr
height2pt&\omit&&\omit&&\omit&&\omit&\cr
\noalign{\hrule}
height2pt&\omit&&\omit&&\omit&&\omit&\cr
& $A$ &&  $-1$ && $1$ &&  $0$ &\cr
height2pt&\omit&&\omit&&\omit&&\omit&\cr
& $B$ && $1$ && $-3$  && $2$  &\cr
height2pt&\omit&&\omit&&\omit&&\omit&\cr
& $D$ && $0 $ && $2$ && $-2$  &\cr
height2pt&\omit&&\omit&&\omit&&\omit&\cr
}  \hrule }
}}
and for the arithmetic genera $p_a(C)$ resp. the intersections $(C,\omega)$ 
with $\omega$ we have: 
\\ \smallskip
\medskip\noindent \vbox{ \bigskip\centerline{\def\quad{\hskip 0.6em\relax}
\def\quod{\hskip 0.5em\relax } \vbox{\offinterlineskip \hrule
\halign{&\vrule#&\strut\quod\hfil#\quad\cr 
height2pt&\omit&&\omit&&\omit&\cr & $C$  &&  $p_a(C)$    &&   $ (C,\omega)
$       &\cr height2pt&\omit&&\omit&&\omit&\cr \noalign{\hrule}
height2pt&\omit&&\omit&&\omit&\cr & $A$ &&  $1$ && $1$ &\cr
height2pt&\omit&&\omit&&\omit&\cr & $B$ && $0	$ && $1$ &\cr
height2pt&\omit&&\omit&&\omit&\cr & $D$ && $0 $ && $0$ &\cr
height2pt&\omit&&\omit&&\omit&\cr }  \hrule } }} 
Assume that $X_k$ is the
special fiber of a semistable hyperelliptic curve $X/S$ with $S$ the spectrum of a
discrete valuation ring $R$ with residue field $k$ and with $X$ regular. Assume
that the Weierstrass points $P_1,\ldots,P_6$ of the generic fiber are rational,
and that their closures in $X$ are distributed over the special fiber as
follows: $P_1,P_2,P_3$ intersect $A$, $P_4$ intersects $B$, and $P_5,P_6$
intersect $D$. The surface $X$ is then the minimal regular model of a genus two
curve given by an equation 
$y^2 = (x-a_1)\cdots(x-a_6)$ 
with $a_1,\ldots,a_6 \in R$ giving
rise to a linear tree $V_0 - V_1 - V_2 - V_3$ with $V_0$ represented by
$a_1,\ldots,a_6$, $V_1,V_2$ represented by $a_4,a_5,a_6$ and $V_3$ represented
by $a_5,a_6$. Thus $\varphi(V_1)=\varphi(V_2)=3$ and $\varphi(V_3)=2$.  The
correspondence with the $P_i$ is via $P_i \leftrightarrow a_i$, the component
$A$ corresponds to $V_0$, the component $B$ corresponds to $V_2$ and the
component $D$ corresponds to $V_3$. The vertex $V_1$ has $C(V)$ equal to $1$. We
compute $e=1$ and Proposition \ref{formulaE} gives $E=A+B+2D$ which is indeed
effective. 

Next, solving the
equation $(2P - \omega + \Phi,C)=0$ for $C=A,B,D$ and demanding that
$(\Phi,P)=0$ one finds  $\Phi_P$ and hence $\Phi_P^2$ for all $P \in
F(S)$. The results are in  the following table: 
\\ \smallskip
\medskip\noindent \vbox{ \bigskip\centerline{\def\quad{\hskip 0.6em\relax}
\def\quod{\hskip 0.5em\relax } \vbox{\offinterlineskip \hrule
\halign{&\vrule#&\strut\quod\hfil#\quad\cr 
height2pt&\omit&&\omit&&\omit&\cr & $P$  &&  $\Phi_P$    &&  
$\Phi_P^2$       &\cr height2pt&\omit&&\omit&&\omit&\cr \noalign{\hrule}
height2pt&\omit&&\omit&&\omit&\cr & $P_1,P_2,P_3$ &&  $-B-D$ && $-1$ &\cr
height2pt&\omit&&\omit&&\omit&\cr & $P_4$ && $-A	$ && $-1$ &\cr
height2pt&\omit&&\omit&&\omit&\cr & $P_5,P_6$ && $-2A-B $ && $-3$ &\cr
height2pt&\omit&&\omit&&\omit&\cr }  \hrule } }} 
The formula in Theorem \ref{main}
reads in our case: 
\[  5 \, d_{X/S} = - \frac{1}{2} \sum_{P \in F(R)} \Phi_P^2 +  9 \,
\delta_{X/S} + 4 \, (E,\omega_{X/S}) \, . \] 
From $E = A+B+2D$ we obtain $(E,\omega_{X/S})=2$. From the tables we further 
read off that  $-\sum_{P \in F(S)}
\Phi_P^2 = 10$. Our description of the special fiber gives 
$\delta_{X/S} = 3$. Finally, one has $d_{X/S} = 8$ by
either \cite[Theorem 3.1]{ka} or \cite[Table \S 5]{ue} (but note that in the
latter reference, the hyperelliptic discriminant is a section of
$\lambda^{\otimes 10}$, not of $\lambda^{\otimes 20}$ as in our set-up, 
so that to get the values
for our $d$, the values for $d$ given in the table in 
\cite{ue} should be multiplied by~$2$). The formula checks.

\section{Weierstrass points on arithmetic surfaces} \label{arithmetic}

In this section we consider Arakelov versions of several of the isomorphisms
considered earlier, associated with Weierstrass points. Let $S$ be a connected 
Dedekind scheme with generic residue characteristic equal to~$0$  and let $\rho
\colon X \to S$ be a   semistable curve of genus $g \geq 1$.  Assume that $X$ is
the regular minimal model of its generic fiber. Under these assumptions, the
Wronskian $Wr$ of $X/S$ is not identically zero and we have a Weierstrass
divisor $W$ flat over $S$ and a residual divisor $E$ supported in the non-smooth
fibers of $X/S$ and in the fibers of positive  characteristic $p$ with $p<2g-1$.
If $\omega$  is the relative dualising sheaf of $\rho$ and $\lambda = \det
\rho_* \omega$ as before, we have a canonical isomorphism $ O_X(W+E) \isom
g(g+1)/2 \cdot \omega - \rho^* \lambda $ of line bundles on $X$ (cf. Proposition
\ref{canonicalisom}). Multiplying both sides by~$4$,  then subtracting
$(g+1)\omega$ from both sides and then pairing with $\omega$ we  obtain from
this a canonical isomorphism $ \langle 4W-(g+1)\omega+4E,\omega \rangle \isom 
\langle (2g-1)(g+1)\omega - 4 \rho^* \lambda, \omega \rangle$. The latter line
bundle was called $Q$ in the proof of Theorem \ref{main}. In that proof, we
applied the Mumford  isomorphism $\langle \omega,\omega \rangle \isom 12 \lambda
- \delta$  to obtain a canonical isomorphism $ Q \isom (3g-1)(8g+4)\lambda -
(2g-1)(g+1)\delta$. Combining both isomorphisms we end up with a canonical
isomorphism: \[ \nu \colon (3g-1)(8g+4)\lambda \isom  \langle
4W-(g+1)\omega+4E,\omega \rangle + (2g-1)(g+1)\delta \] of line bundles on $S$.
Now remark that in the  Faltings-Deligne version of Arakelov theory of
arithmetic surfaces the line bundles under consideration  on $X$ come equipped
with certain canonical hermitian metrics (cf. \cite[Section 2]{fa84}), and
likewise for the line bundles under consideration on $S$  (cf. \cite[Section 3
and 4]{fa84} and \cite[Section 6]{de}).  In particular, if $S=\Spec \cc$  the
isomorphism $\nu$ has a certain norm. Our first result is that this norm is
closely related to an invariant $T$  introduced in \cite{dJ05}. If $\mathcal{X}$
is a compact connected Riemann surface of genus $g \geq 1$, then
$T(\mathcal{X})$ is given by: \begin{align*} T(\mathcal{X}) &= \left( \frac{
\| \vartheta \|(P_1 +\cdots+ P_g-P_{g+1})}{ \prod_{k=1}^g \| \vartheta \|(g P_k
- P_{g+1})^{1/g}} \right)^{2g-2} \cdot \\ &\quad \quad \cdot \left( \frac{
\prod_{k \neq l} \| \vartheta \|(g P_k - P_l)^{1/g}}{ \|J\|(P_1,\ldots,P_g)^2 }
\right) \cdot \prod_{R \in W} \prod_{k=1}^g \| \vartheta \|(g P_k -
R)^{(g-1)/g^4} \, , \end{align*}  where $P_1,\ldots,P_{g+1}$ are generic
points on $\mathcal{X}$, $W$ is the divisor $\div Wr$ of Weierstrass points on
$X$, $\|\vartheta\|$ is the normalised theta function on
$\mathrm{Pic}^{g-1}(\mathcal{X})$ of \cite[p.~401]{fa84}, and $\|J\|$ is a
normalised jacobian determinant involving theta functions  on $\mathrm{Sym}^g
(\mathcal{X})$ introduced in \cite{gu}.  
\begin{prop} \label{normnu} The norm of
$\nu$ for $X$ a smooth proper curve of genus $g \geq 1$ over $S=\Spec \cc$ is
equal to:  \[ (2\pi)^{-4g(2g-1)(g+1)} \cdot T(X(\cc))^{8g^2} \, . \] 
\end{prop}
\begin{proof} We reconstruct the isomorphism $\nu$ again from the isomorphism
\[ O_X(W) \isom g(g+1)/2 \cdot \omega - \rho^* \lambda \] and the Mumford
isomorphism $\langle \omega,\omega \rangle \isom 12 \lambda - \delta$, keeping
track of the norms at each step (note that the residual divisor $E$ now
vanishes).  According to \cite[Th\'eor\`eme 2.2]{mb2} the norm  of the Mumford
isomorphism is equal to $(2\pi)^{4g}e^{-\delta(X(\cc))}$ where $\delta(X(\cc))$
is the Faltings delta-invariant of $X(\cc)$ \cite[p.~402]{fa84}.  The inverse
norm of the isomorphism $O_X(W) \isom g(g+1)/2 \cdot  \omega - \rho^* \lambda$
was called $R(X(\cc))$ in \cite{dJ05} (cf. Definition 5.3 of that paper), and it
follows from \cite[Theorem 4.4]{dJ05} and its proof that $R$, $T$ and $\delta$
are related via: \[ R(X(\cc))^{g-1} = T(X(\cc))^{g^2} \cdot 
e^{-(2g-1)(g+1)\delta(X(\cc))/8}  \] (the reader is warned that the formula in
Theorem~4.4 of \emph{loc. cit.}  contains a misprint: the $g^3$ occurring in the
exponent should be $g^2$). From the isomorphism $O_X(W) \isom g(g+1)/2 \cdot
\omega - \rho^* \lambda$ we get by multiplying both sides by~$4$, by subtracting
$(g+1)\omega$ from both sides and by pairing with $\omega$ the isomorphism \[
\langle 4W-(g+1)\omega,\omega \rangle \isom  \langle (2g-1)(g+1)\omega - 4
\rho^* \lambda, \omega \rangle \] which thus has norm $R(X(\cc))^{-(8g-8)}$.
Composing through with the isomorphism: \[ \langle (2g-1)(g+1)\omega - 4 \rho^*
\lambda, \omega \rangle \isom (3g-1)(8g+4)\lambda - (2g-1)(g+1)\delta \]  which
has norm  $\left((2\pi)^{4g}e^{-\delta(X(\cc))}\right)^{(2g-1)(g+1)}$ we find
the isomorphism: \[ \langle 4W-(g+1)\omega,\omega \rangle \isom 
(3g-1)(8g+4)\lambda - (2g-1)(g+1)\delta \]  having norm: \[ R(X(\cc))^{-(8g-8)}
\cdot  \left((2\pi)^{4g}e^{-\delta(X(\cc))}\right)^{(2g-1)(g+1)} =
(2\pi)^{4g(2g-1)(g+1)} \cdot T(X(\cc))^{-8g^2} \, . \] The norm of $\nu$ is the
inverse of this. 
\end{proof} 
\begin{exam} Let $X$ be a hyperelliptic curve of
genus $g \geq 2$ over $S=\Spec \cc$. It follows from the  proof of Theorem
\ref{main} that the square of the line bundle \[Q=\langle 4W-(g+1)\omega,\omega
\rangle \] on $S$ is canonically isomorphic to the trivial line bundle $O_S$. By
going through a sequence of arguments analogous to those in the proof of
Proposition~2 of \cite{bmmb} it can be verified that the isomorphism $Q^{\otimes
2} \isom O_S$ is in fact an isometry. This implies that the canonical section
$\xi$ of $Q^{\otimes 2}$ has unit norm. According to  \cite[Section 3]{lock} and
\cite[p. 11]{dJ07} the norm of $\Lambda$ satisfies $ \| \Lambda \|^n =
(2\pi)^{4g^2r} \|\Delta_g\|(X(\cc))^g$ where $n={2g \choose g+1}$ and $r={2g+1
\choose g+1}$ and where $\|\Delta_g\|(X(\cc))$ is  the Petersson norm of the
modular discriminant of $X$.  With Proposition \ref{normnu} we obtain the
formula  $T(X(\cc)) = (2\pi)^{-2g} \cdot
\|\Delta_g\|(X(\cc))^{-\frac{3g-1}{8ng}}$ for the $T$-invariant of $X(\cc)$.
\end{exam} 
Proposition \ref{normnu} can be applied to give an explicit formula
for the stable  Faltings height $h(X_K)$ of a smooth proper geometrically
connected  curve $X_K$ of genus $g \geq 2$ over a number field $K$ that has
semistable reduction over $K$. Let $X$ be the regular minimal model of $X_K$
over $S = \Spec O_K$ and let $\omega=\omega_{X/S}$ be the  relative dualising
sheaf of $X/S$. For $P$ a section of $X/S$ denote as before by $\Phi_P$ the
vertical $\qq$-divisor on $X$ such that  $((2g-2)P-\omega+\Phi_P,\Gamma)=0$ for
all irreducible components $\Gamma$ of the geometric fibers of $X/S$ and such
that $P^* \Phi_P$ is trivial. Denote by $h(P)$ the N\'eron-Tate height of the
divisor class of $(2g-2)P-\omega$ in the jacobian of $X_K$. The stable Faltings
height of $X_K$ is defined via $[K\colon \qq]h_F(X_K)= \widehat{\deg} \, \lambda
$, the latter being the  Arakelov degree of $\lambda =\det R\rho_* \omega $
equipped with its canonical hermitian structure at the complex embeddings of $K$
(cf. \cite[Section 3]{fa83}). 
\begin{thm} \label{formulaheight} Assume that the
Weierstrass points on the generic fiber of $X/S$ are rational over $K$. For
$\wp$ a closed point of $S$ let $\delta_\wp$ be the number of singular points in
the geometric fiber at $\wp$. The formula: \begin{align*} (3g-1)(8g+4)
&[K\colon \qq]h_F(X_K)  =  \frac{2[K:\qq]}{g(g-1)}  \sum_{P \in W(S)} h(P) -
\frac{1}{g(g-1)} \sum_{P \in W(S)} \Phi_P^2 \\ & + (2g-1)(g+1) \sum_{\wp \in
|S|} \delta_\wp \log N \wp + 4 (E,\omega_{X/S}) \\ & -    4g(2g-1)(g+1) [K :
\qq] \log (2 \pi) + 8g^2 \sum_{\sigma \colon K \to \cc} \log T(X_\sigma)
\end{align*}  holds. Here for each complex embedding $\sigma$ of $K$  we
denote by $X_\sigma$ the Riemann surface corresponding to the complex curve $X
\times_\sigma \cc$. The intersection numbers $\Phi_P^2$ and $(E,\omega_{X/S})$
should be taken in the Arakelov sense, that is, their local contributions at
each closed point $\wp$ of $S$ should be counted with weight $\log N\wp$. In the
summations over $P \in W(S)$ the Weierstrass points should be counted with their
multiplicity in $W$. 
\end{thm} 
\begin{proof} We start with the canonical
isomorphism: \[ \nu \colon (3g-1)(8g+4)\lambda \isom  \langle
4W-(g+1)\omega+4E,\omega \rangle + (2g-1)(g+1)\delta \] of line bundles on $S$
described at the beginning of this section. The norm of this isomorphism is
provided by Proposition \ref{normnu}. Taking Arakelov degrees on left and right
hand side we then find:  \begin{align*} (3g-1)(8g+4) & \, \widehat{\deg} \,
\lambda =  4(W,\omega) - (g+1)(\omega,\omega) \\ & + (2g-1)(g+1) \sum_{\wp \in
|S|} \delta_\wp \log N \wp + 4 (E,\omega) \\ & -    4g(2g-1)(g+1) [K : \qq] \log
(2 \pi) + 8g^2 \sum_{\sigma \colon K \to \cc} \log T(X_\sigma)  \end{align*} 
where the intersection numbers are Arakelov intersection numbers \cite[Section
2]{fa84}. According to \cite[Section 1.1]{sz} we can write, for any section $P$
of $X/S$: \[ -2\,[K \colon \qq]h(P) = -4g(g-1)(P,\omega) + (\omega,\omega) -
\Phi_P^2 \, . \] Thus, taking the sum over $P$ running through $W(S)$ and
dividing by $g(g-1)$  we find the equality: \[ 4(W,\omega) -
(g+1)(\omega,\omega) = \frac{2[K:\qq]}{g(g-1)} \sum_{P \in W(S)} h(P) -
\frac{1}{g(g-1)} \sum_{P \in W(S)} \Phi_P^2 \, . \] The required formula
follows.  
\end{proof} 
Upon decomposing the terms $h(P)$, $\Phi_P^2$ and
$(E,\omega)$ as sums of local contributions one may use the formula in Theorem
\ref{formulaheight}  to arrive at a provisory definition of a `valuation of the
discriminant' $d_\wp$  for $X/S$ at each closed point $\wp$ of $S$. It might be
worthwhile to study this `valuation of the  discriminant' further especially in
the light of  the fact that, as $\lambda$ is ample on the moduli stack
$\mathcal{M}_g(\cc)$  of complex curves of genus $g \geq 3$, no positive tensor
power $\lambda^{\otimes N}$ of $\lambda$  has a nowhere vanishing global section
on $\mathcal{M}_g(\cc)$.  Specialising $X/S$ to a hyperelliptic curve of genus
$g$, one obtains the valuation of the hyperelliptic discriminant (up to the
factor $3g-1$) as follows from Theorem \ref{main}. An interesting discussion of
the line bundle $\lambda^{\otimes 8g+4}$ on $\mathcal{M}_g(\cc)$ can be found in \cite{hr}.

Most of the terms in our formula for $h_F(X_K)$ can be easily bounded from
below. Indeed, if $P$ is any section of $X/S$ then $h(P) \geq 0$ and $-\Phi_P^2
\geq 0$. Also $(E,\omega)\geq 0$ since $E$ is effective and vertical and the
fibers of $X/S$ do not contain any exceptional curves. By passing, for an
arbitary curve $X_K$ over a number field $K$, to a finite extension $L$ of $K$
such that $X_K$ acquires semistable reduction over $L$ and has all its
Weierstrass points rational over $L$,  we obtain the following corollary.
\begin{cor} \label{lowerboundheight} Let $X_K$ be a smooth proper geometrically
connected  curve of genus $g \geq 2$ over a number field $K$. Then its Faltings
stable height $h_F(X_K)$ satisfies the inequality: \begin{align*} 
(3g-1)&(8g+4)   h_F(X_K)    \geq    - 4g(2g-1)(g+1)  \log (2 \pi) \\  &  + 
(2g-1)(g+1) \frac{1}{[K:\qq]} \sum_{\wp \in |S|} \delta_\wp \log N \wp   + 8g^2
\frac{1}{[K \colon \qq]} \sum_{\sigma \colon K \to \cc} \log  T(X_\sigma) \, .
\end{align*}  Here again for each complex embedding $\sigma$ of $K$ we denote
by  $X_\sigma$ the Riemann surface corresponding to the complex curve $X_K
\times_\sigma \cc$.  
\end{cor}  
It is perhaps interesting to point out a
similarity with an inequality due to Bost \cite[Theorem IV]{bo2}, saying that
the lower bound: \[ (8g+4) h_F(X_K) \geq g \cdot  \frac{1}{[K:\qq]} \left(
\sum_{\wp \in |S|} \delta_\wp \log N \wp + \sum_{\sigma \colon K \to \cc}
\psi(X_\sigma) \right) \] holds, with a function $\psi \colon M_g(\cc) \to \rr$
which is continuous and which satisfies the following logarithmic bound around
the boundary. If $X_0$ is a stable complex curve, if $\rho \colon X \to \Omega$
is its universal deformation with $\Omega$ an open neighbourhood of~$0$ in
$\cc^{3g-3}$,  and if $\delta$ is a local equation for the reduced normal
crossings divisor $\Delta$ in $\Omega$ given by the singular fibers of $\rho$,
then $\psi(X_t) \geq - \log |\delta(t)| + o(\log|\delta(t)|) $ as $t \in \Omega
\setminus \Delta$ goes to~$0$. The `slope' $g/(8g+4)$ featuring in this result 
is much better than the `slope' $(2g-1)(g+1)/(3g-1)(8g+4)$ of our lower bound. 
On the other hand the analytic contributions in our lower bound have the
advantage of being more readily calculable in concrete cases, using only
evaluations of theta functions and their derivatives. An asymptotic analysis of
the function $\log T$ near the boundary of $M_g(\cc)$ was carried out in the
author's thesis, showing that, similar to Bost's result, $T$ can be seen as the
inverse of a distance to the boundary.

An application of the Noether formula for arithmetic surfaces \cite[Th\'eor\`eme
2.5]{mb2} yields a lower bound for the self-intersection of the relative
dualising sheaf.  
\begin{cor} Let $X_K$ be a smooth proper geometrically
connected  curve of genus $g \geq 2$  over a number field $K$. Let $e(X_K)$ be
the normalised self-intersection of the relative dualising sheaf of $X_K$. Then
the inequality: \begin{align*} (3g-1)&(8g+4) e(X_K) \geq  -48g(2g-1)(g+1) 
\log(2\pi) + \frac{8g-8}{[K:\qq]} \sum_{s \in |S|} \delta_\wp \log N \wp \\ &  +
\frac{1}{[K:\qq]} \sum_{\sigma \colon K \to \cc}  \left( 96g^2 \log T(X_\sigma)
- (3g-1)(8g+4)\delta(X_\sigma) \right) \end{align*}  holds.   
\end{cor}
Trying to obtain explicit lower bounds for $e(X_K)$ is interesting in the light
of an effective version of the Bogomolov conjecture \cite{sz90}. Lower bounds of
a similar type appear in \cite[Section 3.3.2]{bu}, which uses more general
Weierstrass points. In a recent preprint \cite{zh} Zhang has obtained explicit
lower bounds for $e(X_K)$  assuming that a certain conjecture of
Grothendieck-Gillet-Soul\'e on the non-negativity of the height of Gross-Schoen
cycles is true. \\

\noindent \textbf{Acknowledgements} The author thanks the Mittag-Leffler
Institute in Djursholm for its hospitality during a visit. He thanks  Sylvain
Maugeais for several discussions related to the theme of this paper, Lidia
Stoppino for pointing out the reference \cite{bo2}, and the anonymous referee
for a number of helpful remarks.  The author is supported by VENI-grant
639.031.619 of  the Netherlands Organisation for Scientific Research (NWO).

\vspace{1cm}

\noindent Address of the author: \\ \\
Mathematical Institute \\
University of Leiden \\
PO Box 9512 \\
2300 RA Leiden \\
The Netherlands \\  
E-mail:  \verb+rdejong@math.leidenuniv.nl+

\end{document}